\documentclass[12pt]{amsart}

\def\frk{\mathfrak}               

\def\mm{{\frk m}}

\def\nn{{\frk n}}

\def\opn#1#2{\def#1{\operatorname{#2}}} 
\opn\chara{char}
\opn\length{\ell}
\opn\pd{pd}
\opn\rk{rk}
\opn\projdim{proj\,dim}
\opn\injdim{inj\,dim}
\opn\rank{rank}
\opn\depth{depth}
\opn\grade{grade}
\opn\height{height}
\opn\embdim{emb\,dim}
\opn\codim{codim}

\opn\Tr{Tr}
\opn\bigrank{big\,rank}
\opn\superheight{superheight}\opn\lcm{lcm}
\opn\trdeg{tr\,deg}%
\opn\reg{reg}
\opn\lreg{lreg}
\opn\div{div}
\opn\Div{Div}
\opn\cl{cl}
\opn\Cl{Cl}
\opn\Spec{Spec}
\opn\Supp{Supp}
\opn\supp{supp}
\opn\Sing{Sing}
\opn\Ass{Ass}
\opn\Ann{Ann}
\opn\Rad{Rad}
\opn\Soc{Soc}
\opn\Ker{Ker}
\opn\Coker{Coker}
\opn\Im{Im}
\opn\Hom{Hom}
\opn\Tor{Tor}
\opn\Ext{Ext}
\opn\End{End}
\opn\Aut{Aut}
\opn\id{id}

\opn\nat{nat}
\opn\pff{pf}
\opn\Pf{Pf}
\opn\GL{GL}
\opn\SL{SL}
\opn\mod{mod}
\opn\ord{ord}
\opn\aff{aff}
\opn\con{conv}
\opn\relint{relint}
\opn\st{st}
\opn\lk{lk}
\opn\cn{cn}
\opn\core{core}
\opn\vol{vol}
\opn\link{link}
\opn\star{star}
\opn\gr{gr}

\def\pot#1#2{#1[\kern-0.28ex[#2]\kern-0.28ex]}

\opn\dirlim{\underrightarrow{\lim}}
\opn\inivlim{\underleftarrow{\lim}}
\let\tensor=\otimes
\let\iso=\cong
\let\Union=\bigcup

\let\Dirsum=\bigoplus

\let\mcone= * 
\let\to=\rightarrow
\let\To=\longrightarrow
\def\Implies{\ifmmode\Longrightarrow \else
     \unskip${}\Longrightarrow{}$\ignorespaces\fi}
\def\implies{\ifmmode\Rightarrow \else
     \unskip${}\Rightarrow{}$\ignorespaces\fi}
\def\iff{\ifmmode\Longleftrightarrow \else
     \unskip${}\Longleftrightarrow{}$\ignorespaces\fi}

\let\:=\colon
\newtheorem{Theorem}{Theorem}[section]
\newtheorem{Lemma}[Theorem]{Lemma}
\newtheorem{Corollary}[Theorem]{Corollary}
\newtheorem{Proposition}[Theorem]{Proposition}

\newtheorem{Definition}[Theorem]{Definition}

\let\epsilon\varepsilon
\let\phi=\varphi
\let\kappa=\varkappa
\opn\Gin{Gin}
\opn\Hilb{Hilb}
\opn\ini{in}
\opn\inim{inm}
\opn\set{set}
\def\pnt{{\raise0.5mm\hbox{\large\bf.}}}
\def\lpnt{{\hbox{\large\bf.}}}

\begin{document}

\title{Sequentially Cohen-Macaulay modules and local cohomology}
\author{J\"urgen Herzog and Enrico Sbarra}
\address{J\"urgen Herzog, Fachbereich Mathematik und Informatik,  
Universit\"at-GHS Essen, 45117 Essen, Germany}
\email{juergen.herzog@@uni-essen.de}
\address{Enrico Sbarra, Fachbereich Mathematik und Informatik,  
Universit\"at-GHS Essen, 45117 Essen, Germany}
\email{enrico.sbarra@@uni-essen.de}

\maketitle

\section*{Introduction}
Let $I\subset R$ be a graded ideal in the polynomial ring $R=K[x_1,\ldots,x_n]$ where $K$ is a 
field, and fix a term order $<$. It has been shown in \cite{Sb} that the Hilbert functions of 
the  local cohomology modules of $R/I$ are bounded by those of $R/\ini(I)$, where $\ini(I)$ 
denotes the initial ideal of $I$ with respect to $<$. In this note we study the question when 
the local cohomology modules of $R/I$ and $R/\ini(I)$ have the same Hilbert function. A 
complete answer to this question can be given for the generic initial ideal $\Gin(I)$ of $I$, 
where $\Gin(I)$ is taken with respect to the reverse lexicographical order and where we assume 
that $\chara(K)=0$. In this case our main result 
(Theorem \ref{main}) says that the local cohomology modules of $R/I$ and $R/\Gin(I)$ have the 
same Hilbert functions if and only if $R/I$ is sequentially Cohen-Macaulay.  

In Section 1 we give the definition of sequentially CM-modules which is due to Stanley 
\cite{St}, and in Theorem \ref{peskine} we present  Peskine's  characterization of 
sequentially CM-modules in terms of Ext-groups. This characterization is used to derive a few 
basic properties of sequentially CM-modules which are needed for the proof of the main result. 

In the following Section 2 we recall some well-known facts about generic initial modules, and 
also prove that $R/\Gin(I)$ is sequentially CM, see Theorem \ref{seqgin}. Section 3 is devoted 
to the proof of the main theorem, and in the final Section 4 we state and prove a squarefree 
version (Theorem \ref{squarefree}) of the main theorem. Its proof is completely different from 
that of the main theorem in the graded case. It is  based upon  a result on componentwise 
linear ideals shown in \cite{AHH1} and the fact (see \cite{HRW}) that the Alexander dual of a 
squarefree componentwise linear ideal defines a sequentially CM simplicial complex.

\section{Sequentially Cohen-Macaulay modules}

We introduce sequentially Cohen-Macaulay modules and  derive some of their basic properties. 
Throughout this section we assume that $R$ is a standard graded Cohen-Macaulay $K$-algebra of 
dimension $n$ with canonical module $\omega_R$.  

The following definition is due to Stanley \cite[Section II, 3.9]{St}.

\begin{Definition}
\label{stanley}
{\em Let $M$ be a finitely generated graded $R$-module. The module $M$ is {\em sequentially 
Cohen-Macaulay} if there exists a finite filtration 
\[
0=M_0\subset M_1\subset M_2\subset\ldots\subset M_r=M
\]
of $M$ by graded submodules of $M$ such that each quotient $M_i/M_{i-1}$ is CM, and  $\dim 
(M_1/M_0)<\dim(M_2/M_1)<\ldots <\dim(M_r/M_{r-1})$.}
\end{Definition}

The following observation follows immediately from the definition:

\begin{Lemma}
\label{modulo}
{\em (a)} Suppose that $M$ is sequentially CM with filtration $0=M_0\subset M_1\subset\ldots 
M_r=M$. Then for any $i=0,\ldots,r$, the module $M/M_i$ is sequentially CM with filtration
$0=M_i/M_i\subset M_{i+1}/M_i\subset\cdots\subset  M_r/M_i$.

{\em (b)} Suppose that $M_1\subset M$ and $M_1$ is CM, and $M/M_1$ is sequentially CM with 
$\dim M_1<\dim M/M_1$. Then $M$ is sequentially CM.
\end{Lemma}

In order to simplify notation we will write $E^i(M)$ for $\Ext^i_R(M,\omega_R)$.

\begin{Proposition}
\label{exti}
Suppose that $M$ is sequentially CM with a filtration as in {\em \ref{stanley}}, and assume 
that $d_i=\dim M_i/M_{i-1}$. Then 
\begin{enumerate}
\item[(a)]  $E^{n-d_i}(M)\iso E^{n-d_i}(M_i/M_{i-1})$,  and  
is CM of dimension $d_i$ for $i=1,\ldots,r$, and 
$E^j(M)=0$ if $j\not\in\{n-d_1,\ldots,n-d_r\}$.
\item[(b)] $E^{n-d_i}(E^{n-d_i}(M))\iso M_i/M_{i-1}$ for $i=1,\ldots,r$.
\end{enumerate}
\end{Proposition}  

\begin{proof}
(a) We proceed by induction on $r$. From the short exact sequence $0\to M_1\to M\to M/M_1\to 
0$ we obtain the long exact sequence
\[
\cdots\To E^j(M/M_1)\To E^j(M)\To E^j(M_1)\to E^{j+1}(M/M_1)\To\cdots
\]
From \cite[Theorem 3.3.10]{BH} it follows that $E^j(M_1)=0$ if $j\neq n-d_1$, and that 
$E^{n-d_1}(M_1)$ is CM of dimension $d_1$. Thus we get an exact sequence
\begin{eqnarray}
\label{semishort}
0&\to E^{n-d_1}(M/M_1)\to E^{n-d_1}(M)\to E^{n-d_1}(M_1)\\
&\to E^{n-d_1+1}(M/M_1)\to E^{n-d_1+1}(M)\to 0,\nonumber
\end{eqnarray}
and isomorphisms $E^j(M/M_1)\iso E^j(M)$ for all $j\neq n-d_1,n-d_1+1$.

By Lemma \ref{modulo}, the module $M/M_1$ is sequentially CM and has a CM filtration of length 
$r-1$. Hence by induction hypothesis we have $E^{n-j}(M/M_1)=0$ for $j\neq 
\{d_2,\ldots,d_r\}$.  This implies  that $E^{n-d_1}(M/M_1)=E^{n-d_1+1}(M/M_1)=0$, and hence by 
(\ref{semishort}) we have $E^{n-d_1}(M)\iso E^{n-d_1}(M_1)$, and $E^{n-d_1+1}(M)=0$.
Summing up we conclude that $E^{n-d_1}(M)\iso E^{n-d_1}(M_1)$ and $E^j(M)\iso E^j(M/M_1)$ for 
$j\neq n-d_1$. Thus the assertion follows from the induction hypothesis and the fact that 
$E^{n-d_1}(M_1)$ is CM of dimension $d_1$.

(b) follows from (a) and \cite[Theorem 3.3.10]{BH} since for any CM-module $N$ of dimension 
$d$ one has 
$N\iso E^{n-d}(E^{n-d}(N))$.
\end{proof}

It is quite surprising that \ref{exti} has a strong converse. The following theorem is due to 
Peskine. Since there is no published proof available we present here a proof for  the 
convenience of the reader.

\begin{Theorem}
\label{peskine}
The following two conditions are equivalent:
\begin{enumerate}
\item[(a)] $M$ is sequentially CM;
\item[(b)] for all $0\leq i\leq \dim M$, the modules $E^{n-i}(M)$ are either $0$ or CM of 
dimension $i$.
\end{enumerate}
\end{Theorem}

The implication (a)\implies (b) follows from \ref{exti}. For the other direction we first need 
to show
\begin{Lemma}
\label{crucial}
Let $t=\depth M$, and suppose that $E^{n-t}(M)$ is CM of dimension $t$. Then there exists a 
natural monomorphism $\alpha: E^{n-t}(E^{n-t}(M))\to M$, and the induced map 
$E^{n-t}(\alpha)\: E^{n-t}(M)\to E^{n-t}(E^{n-t}(E^{n-t}(M)))=E^{n-t}(M)$ is an isomorphism.
\end{Lemma}

\begin{proof}
We write $R=S/I$, where $S$ is a polynomial ring. Let $\mm$ be the graded maximal ideal of 
$R$, and $\nn$ the graded maximal ideal of $S$. By the Local Duality Theorem (see 
\cite[Theorem 3.6.10]{BH}) we have
\begin{eqnarray}
\label{local1}
\Ext^i_R(M,\omega_R)\iso \Hom_R(H^{n-i}_\mm(M),E_R(K)), 
\end{eqnarray}
and 
\begin{eqnarray}
\label{local2}
\Ext^i_S(M,\omega_S)\iso \Hom_S(H^{m-i}_\nn(M),E_S(K)), 
\end{eqnarray} 
where $m=\dim S$. Since $H^i_{\nn}(M)\iso H^i_\mm(M)$, and since $\Hom_S(R,E_S(K))\iso 
E_R(K)$, we see that 
\begin{eqnarray*}
\Hom_S(H^j_\nn(M),E_S(K))&\iso& \Hom_S(H^j_\mm(M),E_S(K))\\
&\iso &\Hom_R(H^j_\mm(M),\Hom_S(R, E_S(K))\\
&\iso& \Hom_R(H^{j}_\mm(M),E_R(K)).
\end{eqnarray*}
Therefore (\ref{local1}) and (\ref{local2}) imply that
\[
\Ext^{n-t}_R(M,\omega_R)\iso\Ext^{m-t}_S(M, \omega_S),
\]
and we hence may as well assume that $R$ is a polynomial ring.

Let 
\[
F_\lpnt\: 0\To F_{n-t}\To \cdots\To F_1\To F_0 \To 0
\]
be the minimal graded free resolution of $M$. Note that $\omega_R=R(-n)$ since $R$ is a 
polynomial ring.  Then the $\omega_R$-dual $F^*_\lpnt$ of $F_\lpnt$ is the free complex
\[
0\To F^*_{n-t}\To \cdots \To F^*_1\To F^*_0\To 0 
\]
with $F^*_i=\Hom_R(F_{n-i},\omega_R)$, and $H_0(F^*)=E^{n-t}(M)$.

Let $G_\lpnt$ be the minimal graded free resolution of $E^{n-t}(M)$.  Then there exists a 
comparison map $\phi_\lpnt\: F^*_\lpnt\To G_\lpnt$ which extends the identity on 
$H_0(F^*_\lpnt)=E^{n-t}(M)=H_0(G_\lpnt)$.

Since by assumption $E^{n-t}(M)$ is CM of dimension $t$, the complex $G_\lpnt$ has the same 
length  as $F^*_\lpnt$, namely $n-t$. Thus the $\omega_R$-dual $\phi^*_\lpnt\: G^*_\lpnt\to 
F_\lpnt$ of $\phi_\lpnt$ induces a natural homomorphism $\alpha=H_0(\phi^*_\lpnt)\: 
E^{n-t}(E^{n-t}(M))=H_0(G^*_\lpnt)\to H_0(F_\lpnt)=M$. Here $G^*_i=\Hom_R(G_{n-i},\omega_R)$ 
and $\phi^*_i=\Hom_R(\phi_{n-i},\omega_R)$ for all $i$.

Since $E^{n-t}(M)$ is CM by assumption, the complex $G^*_\lpnt$ is exact, and hence a free 
resolution of $E^{n-t}(E^{n-t}(M))$, and so the induced map $E^{n-t}(E^{n-t}(E^{n-t}(M)))\to 
E^{n-t}(M)$ is given by $H_0(\phi^{**})=H_0(\phi)=\id$.

It remains to be shown that $\alpha$ is a monomorphism.  Let $C_\lpnt$ be the mapping cone of 
$\phi^*_\lpnt:\ G^*_\lpnt\to F_\lpnt$. Since $F_\lpnt$ and $G^*_\lpnt$ are acyclic, it follows 
that $H_1(C_\lpnt)\iso \Ker(\alpha)$ and $H_i(C_\lpnt)=0$ for $i>1$. Notice that 
$\phi^*_{n-t}$ is an isomorphism, since this is the case for $\phi_0$. Hence the chain map 
$C_{n-t+1}\to C_{n-t}$ is split injective, and so by cancellation we get 
a new complex of free $R$-modules 
\[
\tilde{C}_\lpnt: 0\To D_{n-t}\to C_{n-t-1}\To\cdots\To C_0\To 0,
\]
where $D_{n-t}=\Coker(C_{n-t+1}\to C_{n-t})$. Again we have $H_1(\tilde{C}_\lpnt)\iso 
\Ker(\alpha)$ and $H_i(\tilde{C}_\lpnt)=0$ for $i>1$.

Now suppose that $\Ker(\alpha)\neq 0$, and let $P$ be a minimal prime ideal of the support of 
$\Ker(\alpha)$. Since $\Ker(\alpha)\subset E^{n-t}(E^{n-t}(M))$, and since 
$E^{n-t}(E^{n-t}(M))$ is a CM-module of dimension $t$, it follows that $P$ is a minimal prime 
ideal of $E^{n-t}(E^{n-t}(M))$ with  $\height P=n-t$. Therefore 
$L_\lpnt=\tilde{C}_\lpnt\tensor R_P$, is a complex of length $n-t$ with $\depth(L_i)=n-t$ for 
all $i$,  $\depth(H_1(L_\lpnt))=0$ and $H_i(L_\lpnt)=0$ for $i>0$. By the Peskine-Szpiro lemme 
d'acyclicit\'e \cite{PS} this implies that $\tilde{C}$ is acyclic, a contradiction.  
\end{proof}

\begin{proof}[Proof of \ref{peskine}] We proceed by induction on $n-t$. 
Let $t=\depth M$, and let $M_1$ be the image of $E^{n-t}(E^{n-t}(M))\to M$. By \ref{crucial}, 
the module $M_1$ is a CM-module of dimension $t$. 

Consider the short exact sequence $0\to M_1\to M\to M/M_1\to 0$. As in the proof of \ref{exti} 
we get an exact sequence 
\begin{eqnarray*}
0&\to E^{n-t}(M/M_1)\to E^{n-t}(M)\to E^{n-t}(M_1)\\
&\to E^{n-t+1}(M/M_1)\to E^{n-t+1}(M)\to 0,\nonumber
\end{eqnarray*}
and isomorphisms $E^j(M/M_1)\iso E^j(M)$ for all $j\neq n-t,n-t+1$.\\
\noindent
Since $E^{n-t}(M)\to E^{n-t}(M_1)$ is an isomorphism  (see \ref{crucial}), we deduce from the 
above exact sequence that $E^{n-t}(M/M_1)=0$, and that  $E^{n-t+1}(M/M_1)\iso E^{n-t+1}(M)$. 
Thus we have $E^j(M/M_1)\iso E^j(M)$ for $j<n-t$, and $E^j(M/M_1)=0$ for $j\geq n-t$. Hence, 
by induction hypothesis, $M/M_1$ is sequentially CM, and so is $M$ by \ref{modulo}.
\end{proof}

An immediate application of \ref{peskine} is

\begin{Corollary}
\label{direct}
Let $M$ be a finite direct sum of sequentially CM-modules. Then $M$ is sequentially $CM$.
\end{Corollary}
   
As a consequence of \ref{exti} and \ref{crucial} we get 

\begin{Corollary}
\label{unique}
A filtration of a sequentially CM module satisfying the conditions of {\em \ref{stanley}} is 
uniquely determined.
\end{Corollary}

\begin{proof} 
Let $t=\depth M$. The first module $M_1$ in the filtration must be the image of 
$E^{n-t}(E^{n-t}(M))\to M$. Then one makes use of an induction argument to $M/M_1$ to obtain 
the desired result.
\end{proof}

Notice that $M_1=H^0_\mm(M)$ if $\depth M=0$. Thus, \ref{stanley} together with \ref{modulo} 
imply

\begin{Corollary}
\label{depth0}
An $R$-module $M$ is sequentially CM if and only if $M/H^0_\mm(M)$ is sequentially CM. 
\end{Corollary}

In what follows we denote by  $E^\pnt(M)=\Dirsum_iE^i(M)$. Then we get

\begin{Corollary}
\label{x}
Suppose that $x\in R$ is a homogeneous $M$- and $E^\pnt(M)$-regular element. Then $M$ is 
sequentially CM if and only if $M/xM$ is sequentially CM.
\end{Corollary}  

\begin{proof} Since $x$ is $E^\pnt(M)$-regular, the long exact Ext-sequence derived from 
\[
0\To M(-1)\buildrel x\over\To  M\To M/xM\to 0
\]
splits into short exact sequences
\[
0\To E^{n-i}(M)\buildrel x\over \To E^{n-i}(M)(1)\To E^{n-i+1}(M/xM)\To 0.
\]
It follows that $E^{n-i}(M)$ is CM of dimension  $i$ if and only if $E^{n-i+1}(M/xM)$ is CM of 
dimension  $i-1$. Thus  \ref{peskine} implies the assertion.  
\end{proof} 

In conclusion we would like  to remark that the same theory is valid in the category of 
finitely generated $R$-modules, where $R$ is a local  CM ring and a factor ring of a  regular 
local ring.

\section{Generic initial modules}

In this section we recall a few facts on generic initial modules, which are mostly due to 
Bayer and Stillman, and can be found in \cite{Ei}. 

Let $R=K[x_1,\ldots,x_n]$ be the polynomial ring over a  field $K$ of characteristic $0$, and 
let $M$ be a graded module with graded free presentation $F/U$. Throughout this section  let 
$<$ be a term order that refines the partial order by degree and that satisfies 
$x_1>x_2>\cdots >x_n$. We fix a graded basis $e_1,\ldots, e_m$ of $F$,  and extend the order 
$<$ to $F$ as follows: Let $ue_i$ and $ve_j$ be  monomials (i.e.\ $u$ and $v$ are  monomials 
in $R$). We set $ue_i>ve_j$ if either $\deg(ue_i)>\deg(ve_j)$, or the degrees are the same and 
$i<j$, or $i=j$ and $u>v$.

We set 
\[
U^{sat}=\Union_rU:\mm^r=\{f\in F\: f\mm^r\in U\quad\text{for some $r$}\}.
\]

From now on let $<$ denote the reverse lexicographic order. In the next proposition we collect 
all the results which will be needed later.
\begin{Proposition}
\label{needed}
For generic choice of coordinates one has:
\begin{enumerate}
\item[(a)] $\dim F/\Gin(U)=\dim F/U$ and $\depth F/\Gin(U)=\depth F/U$;
\item[(b)] $x_n$ is $F/U$ regular if and only if $x_n$ is $F/\Gin(U)$ regular;
\item[(c)] $\Gin(U)^{sat}=\Gin(U^{sat})$.
\end{enumerate}
\end{Proposition}

\begin{proof}
After  a generic choice of coordinates we may assume that $\Gin(U)=\ini(U)$. The first 
statement in (a) is true for any term order, while the second statement about the depth  and 
assertion (b) follow from \cite[Theorem 15.13]{Ei} because we may assume that the sequence 
$x_n, x_{n-1},\ldots, x_{n-t+1}$ is $M$-regular if $\depth M=t$.

By the module version of \cite[Proposition 15.24]{Ei}, and by \cite[Proposition 15.12]{Ei} one 
has that 
\[
\Gin(U)^{sat}= \Union_r(\Gin(U):x_n^r)=\Union_r\Gin(U:x_n^r).
\] 
On the other hand for a generic choice of coordinates we have $U^{sat}=\Union_r (U:x_n^r)$. 
Therefore $\Union_r\Gin(U:x_n^r)=\Gin(\Union_r (U:x_n^r))=\Gin(U^{sat})$, which yields the 
last assertion.
\end{proof}

The following result will be crucial for the proof of the main theorem of this paper.

\begin{Theorem}
\label{seqgin}
The module $F/\Gin(U)$ is sequentially CM.
\end{Theorem}

\begin{proof}
Observe that, since we assume $\chara(K)=0$, we have $\Gin(U)=\Dirsum_j I_je_j$ where for each 
$j$, $I_j$ is a strongly stable ideal , cf.\ \cite[Theorem 15.23]{Ei}. Hence 
$F/\Gin(U)\iso\Dirsum_jR/I_j$, so that, by \ref{direct} one only has to prove that $R/I$ is 
sequentially CM for any strongly stable ideal $I\subset R$. 

Recall that a monomial ideal is strongly stable if 
for all $u\in G(I)$ and all $i$  such that $x_i$ divides $u$ one has $x_j(u/x_i)\in I$ for all 
$j<i$. Here $G(I)$ denotes the unique minimal set of monomial generators of $I$. 

For a monomial $u$ we let $m(u)=\max\{i\: x_i\ \text{divides}\ u\}$, and $s=\max\{m(u)\:u\in 
G(I)\}$. Let $R'=K[x_1,\ldots,x_s]$, and let $J\subset R'$ be the unique monomial ideal such 
that $I=JR$. It is clear that $J$ is a strongly stable ideal in $R'$. Thus it follows that 
$J^{sat}=  \Union_r(J:x_s^r)$. Note that $J^{sat}$ contains $J$ properly and is  strongly 
stable. Let
$I_1=\Union_r(I:x_s^r)$. Then $I_1=J^{sat}R$, and since the extension $R'\to R$ is flat, we 
have $I_1/I\iso (J^{sat}/J)\tensor_{R'}R$. Now $J^{sat}/J$ is a non-trivial $0$-dimensional CM 
module over $R'$, and therefore $M_1=I_1/I\subset R/I$ is an $(n-s)$-dimensional CM-module 
over $R$. Next we observe that $(R/I)/M_1=R/I_1$ and that $I_1$ is strongly stable. Since 
$\dim R/I_1\geq n-\max\{m(u)\: u\in G(I_1)\}>n-s$, the assertion of the theorem follows from 
\ref{modulo}.
\end{proof}

\section{The main theorem}

As in the previous section we let $K$ be a field of characteristic 0, $R=K[x_1,\ldots,x_n]$ be 
the polynomial ring over $K$  and $M$ be a finitely generated graded $R$-module with graded 
free presentation $M=F/U$. We want to compare the Hilbert functions of the local cohomology 
modules of $F/U$ and $F/\Gin(U)$, where $\Gin(U)$ is taken with respect to the reverse 
lexicographical order. In general one has (see \cite{Sb}) a coefficientwise inequality  
$\Hilb(H^i_\mm(F/U))\leq \Hilb(H^i_\mm(F/\Gin(U)))$. The main purpose of this section is to 
prove

\begin{Theorem}
\label{main}
The following conditions are equivalent:
\begin{enumerate}
\item[(a)] $F/U$ is sequentially CM;
\item[(b)] for all $i\geq 0$ one has $\Hilb(H^i_\mm(F/U))=  \Hilb(H^i_\mm(F/\Gin(U)))$.
\end{enumerate}
\end{Theorem}

\begin{proof} (a)\implies (b): Set $M=F/U$ and $N=F/\Gin(U)$. We proceed by induction on $\dim 
M$. Suppose $\dim M=0$, then $\dim N=0$ and $\Hilb(M)=\Hilb(N)$. Since
$H^0_\mm(M)=M$, $H^0_\mm(N)=N$ and $H^i_\mm(M)=H^i_\mm(N)=0$ for $i>0$, the assertion follows 
in this case.    

Now suppose that $\dim M>0$. Assume first that $\depth M=0$. We have $M/H^0_\mm(M)\iso 
F/U^{sat}$, and, by \ref{needed},  $N/H^0_\mm(N)=F/\Gin(U)^{sat}=F/\Gin(U^{sat})$. By 
\ref{depth0}, we also know that $M/H^0_\mm(M)$ is sequentially CM. Thus, if the implication 
(a)\implies (b) were known for modules of positive depth, it would follow that 
\[
\Hilb(H^i_\mm(M))=\Hilb(H^i_\mm(M/H^0_\mm(M)))=\Hilb(H^i_\mm(N/H^0_\mm(N)))=
\Hilb(H^i_\mm(M))
\]
for all $i>0$. Notice that $H^0_\mm(M)=U^{sat}/U$ and $H^0_\mm(N)=\Gin(U^{sat})/\Gin(U)$. 
However, since $M=F/U$ and $N=F/\Gin(U)$, and since $F/U^{sat}$ and $F/\Gin(U^{sat})$ have the 
same Hilbert function, we conclude that also $H^0_\mm(M)$ and $H^0_\mm(N)$ have the same 
Hilbert function. 

These considerations show that we may assume that $\depth M>0$.
Accordingly, $\depth N>0$ by \ref{needed}, and $N$ is sequentially CM by \ref{seqgin}. Since 
$M$ and $N$ are sequentially CM we have $\depth E^\pnt(M)>0$ and $\depth E^\pnt(N)>0$. We may 
assume  that the coordinates are chosen generically so that $\Gin(U)=\ini(U)$ and that $x_n$ 
regular on $E^\pnt(M)$ and regular on $E^\pnt(N)$. According to \ref{x}, $M/x_nM=F/(U+x_nF)$ 
is  sequentially CM. Therefore our induction hypothesis, together with \cite[Proposition 
15.12]{Ei}, implies that the Hilbert functions of the local cohomology modules of $M/x_nM$ and 
of $F/\Gin(U+x_nF)=F/(\Gin(U)+x_nF)=N/x_nN$ are the same.

We have short exact sequences
of graded $R$-modules
\[
0\To H^{i-1}_\mm(M/x_nM)\To H^i_\mm(M)(-1)\buildrel x_n\over \To H^i_\mm(M)\To 0,
\]
and 
\[
0\To H^{i-1}_\mm(N/x_nN)\To H^i_\mm(N)(-1)\buildrel\ x_n\over \To H^i_\mm(N)\To 0,
\]
because $x_n$ is regular on $E^\pnt(M)$ and $E^\pnt(N)$.
Therefore applying the induction hypothesis to $M/x_nM$ we get 
\begin{eqnarray*}
\Hilb(H^i_\mm(M))(t-1)&=&\Hilb(H^{i-1}_\mm(M/x_nM))\\
&=&\Hilb(H^{i-1}_\mm(N/x_nN))=
\Hilb(H^i_\mm(N))(t-1),
\end{eqnarray*}
from which we deduce that $\Hilb(H^i_\mm(M))=\Hilb(H^i_\mm(N))$.

(b)\implies (a): We proceed again by induction on $\dim M$. With the same arguments as in the 
proof of the first implication, we may assume that $\depth M>0$. Therefore $\depth N>0$, too, 
and since we are working with  generic coordinates and $N$ is sequentially CM by \ref{seqgin}, 
$x_n$ is $E^\pnt(N)$-,  $M$- and $N$-regular. We shall show that $x_n$ is also 
$E^\pnt(M)$-regular. Since $x_n$ is $E^\pnt(N)$-regular, the long exact cohomology sequence 
derived from  $0\to N(-1)\buildrel x_n\over \to  N\to N/x_nN\to 0$ splits into short exact 
sequences
\begin{eqnarray}
\label{short}
0\To H^{i-1}_\mm(N/x_nN)\To H^i_\mm(N)(-1)\buildrel\ x_n\over \To H^i_\mm(N)\To 0,
\end{eqnarray}
We show by induction on $i$ that the corresponding sequences for $M$ are also exact.
For $i=0$ the assertion is trivial, since $H^0_\mm(M)=0$. Now let $i>0$, and assume that the 
assertion is true for all $j<i$. Then $ H^{i-1}_\mm(M)(-1)\buildrel\ x_n\over \to 
H^{i-1}_\mm(M)$ is surjective, and we obtain the exact sequence
\[
0\To H^{i-1}_\mm(M/x_nM)\To H^i_\mm(M)(-1)\buildrel\ x_n\over \To H^i_\mm(M). 
\]
Suppose multiplication with $x_n$ is not surjective, then there exists a degree $a$ such that  
$H^i_\mm(M)_{a-1}\buildrel\ x_n\over \To H^i_\mm(M)_a$ is not surjective. Using  
(\ref{short}), and the hypothesis that the local cohomology modules of $M$ and $N$ have the 
same Hilbert function, one has  
\begin{eqnarray*}
\dim_K  H^{i-1}_\mm(M/x_nM)_a&>& \dim_K H^i_\mm(M)_{a-1}-\dim H^i_\mm(M)_a\\
&=&\dim_K H^i_\mm(N)_{a-1}-\dim H^i_\mm(N)_a\\
&=& \dim_K  H^{i-1}_\mm(N/x_nN)_a.
\end{eqnarray*}
This is a contradiction, since $M/x_nM=F/(U+x_nF)$ and $N/x_nN=F/\Gin(U+x_nF)$, and consequently  
$\dim_KH^{i-1}_\mm(M/x_nM)\leq \dim_KH^{i-1}_\mm(N/x_nN)$. 

Now it follows that $x_n$ is $E^\pnt(M)$-regular, and also that the cohomology modules of 
$M/x_nM$  and $N/x_nN$ have the same Hilbert functions, whence our induction hypothesis 
implies that $M/x_nM$ is sequentially CM. Since $x_n$ is $E^\pnt(M)$-regular, we  finalöly 
deduce that $M$ is sequentially CM.
\end{proof}

\section{The squarefree case}

In this section we will state and prove the squarefree analogue of Theorem \ref{main}.
Let $\Delta$ be a simplicial complex on the vertex set $[n]=\{1,\ldots,n\}$, 
and let $I_\Delta\subset R$ be the Stanley-Reisner ideal of $\Delta$, where 
$R=K[x_1,\ldots,x_n]$ and $K$ is field of characteristic $0$.  The $K$-algebra 
$K[\Delta]=R/I_\Delta$ is the Stanley-Reisner ring of $\Delta$. 

We recall the concept of symmetric algebraic shifting which was introduced by Kalai in  
\cite{Ka2}: Let $u\in R$ be a monomial, $u=x_{i_1}x_{i_2}\cdots x_{i_d}$ with $i_1\leq i_2\leq 
\ldots\leq i_d$. We define 
\[
u^\sigma =x_{i_1}x_{i_2+1}\cdots x_{i_d+d-1}.
\]
Note that $u^\sigma$ is a squarefree monomial (in a possibly bigger polynomial ring). 

As usual  the unique minimal monomial set of generators of a monomial ideal $I$ is denoted by 
$G(I)$. 

The {\em symmetric algebraic shifted complex of $\Delta$} is the simplicial complex $\Delta^s$ 
whose Stanley-Reisner ideal $I_{\Delta^s}$ is generated by the squarefree monomials 
$u^\sigma$ with $u\in \Gin(I_\Delta)$.

We quote the following properties of  $I_{\Delta^s}$ from \cite{AHH} and \cite{AHH1}:
\begin{enumerate}
\item[(i)] $I_{\Delta^s}$ is a strongly stable ideal in $R$;
\item[(ii)] one has the following inequality of graded Betti numbers:
\[
\beta_{ij}(I_\Delta)\leq \beta_{ij}(I_{\Delta^s});
\]
\item[(iii)] $I_\Delta$ and $I_{\Delta^s}$ have the same graded Betti numbers if and only if 
$I_\Delta$ is componentwise linear. 
\end{enumerate}

Recall that an ideal $I\subset R$ is called {\em componentwise linear} if in each degree $i$, 
the ideal generated by the $i$-th graded component $I_i$ of $I$ has a linear resolution. 

Let $\Delta^*$ denote the Alexander dual of $\Delta$, i.e., the simplicial complex 
\[
\Delta^*=\{F\subset [n]\: [n]\setminus F\not\in \Delta\}.
\]
It has been noted in \cite[Theorem 9]{HRW} that
\[
K[\Delta] \quad\text{is sequentially CM}\iff I_{\Delta^*}\quad\text{is componentwise linear}.
\]

\begin{Theorem}
\label{squarefree}
Let $\Delta$ be a simplicial complex. Then 
\begin{enumerate}
\item[(a)] $\Hilb(H^i_\mm(K[\Delta]))\leq \Hilb(H^i_\mm(K[\Delta^s]))$ for all $i$.
\item[(b)] The local cohomology module of $K[\Delta]$ and $K[\Delta^s]$ have the same Hilbert 
function if and only if $K[\Delta]$ is sequentially CM.
\end{enumerate}
\end{Theorem}

\begin{proof}
Part (a) is proved in \cite{Sb}. For the proof of (b) we shall need the following result which 
also can be found in \cite{Sb}: for all $i\geq 0$ and $j\geq 0$ one has
\begin{eqnarray}
\label{enrico}
\dim_K H^i_\mm(K[\Delta])_{-j}=\sum_{h=0}^n\tbinom{n}{h}\tbinom{h+j-1}{j}\beta_{i-h+1, 
n-h}(K[\Delta^*]).
\end{eqnarray}
(Observe that $H^i_\mm(K[\Delta])_j=0$ for $j> 0$ and all $i$, as shown by Hochster, see 
\cite{Ho} and \cite[Theorem 5.3.8]{BH}).

Now suppose that $K[\Delta]$ is sequentially CM. Then $I_{\Delta^*}$ is componentwise linear, 
and hence $\beta_{ij}(K[\Delta^*])=\beta_{ij}(K[(\Delta^*)^s])$ by Property (iii) of symmetric 
algebraic shifting. Since $(\Delta^*)^s=(\Delta^s)^*$, Formula (\ref{enrico}) shows that 
$\dim_K H^i_\mm(K[\Delta])_{-j}=\dim_K H^i_\mm(K[\Delta^s])_{-j}$ for all $i$ and $j$, as 
desired. 

For the viceversa, let $H$ be the $(n+1)\times (n+1)$-matrix with entries 
$h_{ij}=\dim_KH^i_\mm(K[\Delta])_{-j}$, $i,j=0,\ldots,n$, $B$ the $(n+1)\times (n+1)$-matrix 
with entries
$b_{hi}=\tbinom{n}{h}\beta_{i-h+1,n-h}(K[\Delta^*])$ and $A$ the $(n+1)\times (n+1)$-matrix 
with entries
$a_{jh}=\tbinom{h+j-1}{j}$. Then (\ref{enrico}) says that $H^t=AB$. Since $A$ is invertible, we 
see that the numbers $b_{hi}=\tbinom{n}{h}\beta_{i-h+1,n-h}$ are determined by the Hilbert 
functions of the local cohomology modules of $K[\Delta]$. Thus, if    
$\dim_K H^i_\mm(K[\Delta])_j=\dim_K H^i_\mm(K[\Delta^s])_j$ for all $i$ and $j$, then 
the numbers $b_{hi}$ for $\Delta^*$ and $(\Delta^*)^s$ coincide, which in turn implies that 
their graded Betti numbers are the same (because $\beta_{ij}=b_{n-j,i+n-j-1}/\tbinom{n}{n-j}$). 
By Property (iii) of symmetric algebraic shifting this implies that $I_\Delta^*$ is 
componentwise linear, and hence $K[\Delta]$ is sequentially CM.  
\end{proof}


\newpage


\begin{thebibliography}{99}

\bibitem{AHH} A.\ Aramova, J.\ Herzog, and T.\ Hibi, ``Shifting
operations and graded Betti numbers'', to appear in J. of Alg.\ Comb.

\bibitem{AHH1} A.\ Aramova, J.\ Herzog, and T.\ Hibi, ``Ideals
with stable Betti
numbers'', to appear in Adv.\ Math.

\bibitem{AAH} A.\ Aramova, L.L.\ Avramov, J.\ Herzog,
``Resolutions of monomial ideals and cohomology over exterior algebras, to appear in TAMS.

\bibitem{B} D.\ Bayer, ``The division algorithm and the Hilbert
scheme'', Ph.D. Thesis, Harvard 1982. 

\bibitem{BS} D.\ Bayer, and M.\ Stillman, ``A theorem on refining
division orders by the reverse lexicographical order'', {\it Duke
Math. J.}  {\bf 55} (1987),  321 -- 328.


\bibitem{BH} W.\ Bruns, J.\ Herzog, {\em Cohen-Macaulay
rings},  Revised Edition, Cambridge, 1996. 

\bibitem{Ei} D.\ Eisenbud, {\em Commutative algebra, with a view
towards 
algebraic geometry}, Graduate Texts Math., Springer. 1995.

\bibitem{EK} S.\ Eliahou, M.\ Kervaire, ``Minimal resolutions of
some monomial ideals'', {\em J.\ Alg.} {\bf 129} (1990), 1--25.


\bibitem{G} A. Galligo, ``Apropos du th\'eor\`eme de pr\'eparation
de Weierstrass'', {\it in} ``Fonctions de plusiers variables
complexes'',
Springer Lect.\ Notes Math. {\bf 409} (1974), 543 -- 579 .  

\bibitem{Gre} M.\ Green, ``Generic initial ideals'', {\em
Proc.\ CRM -96, 
Six Lectures on Commutative Algebra}, 
Barcelona, Spain, {\bf 166}, Birkh\"auser, 1998, 119 -- 186, . 

\bibitem{HRW} J.\ Herzog, V.\ Reiner and V.\ Welker, ``Componentwise linear ideals and Golod 
rings''


\bibitem{HP} J.\ Herzog and  D.\ Popescu, ``On the regularity of
$p$-Borel ideals'', preprint 1999. 

\bibitem{Ho} M.\ Hochster, ``Cohen--Macaulay rings, combinatorics, 
and simplicial complexes'', {\em Proc.\  Ring Theory II}, Lect. 
Notes in Pure and Appl Math., {\bf 26}, Dekker,
New York, 1977, 171 -- 223. 

\bibitem{Ka1} G.\ Kalai,  ``Algebraic shifting'', Unpublished
manuscript, 1993.

\bibitem{Ka2} G.\ Kalai,  ``The diameter of graphs of convex
polytopes and
$f$-vector theory'', {\em Proc.\ Applied Geometry and
Discrete Mathematics},
DIMACS Series in Discrete Mathematics and Theoretical Computer
Science, {\bf 4}, Amer. Math. Soc., 1991, 387--411.

\bibitem{PS} C.\ Peskine, L.\ Szpiro, ``Dimension projective finite cohomology locale´´,
Publ.\ Math.\ I.H.E.S., {\bf 42} (1972), 47--119.


\bibitem{Sb} E.\ Sbarra, ``Upper bounds for local cohomology for
rings with given Hilbert function'', to appear in Comm.\ in Alg.

\bibitem{St} R.\ P.\ Stanley, {\em Combinatorics and Commutative
Algebra,} Birkh\"{a}user,  1983.

\end{thebibliography}
\end{document}